\documentclass[12pt,a4paper,leqno]{amsart}

\usepackage[latin1]{inputenc}
\usepackage[T1]{fontenc}
\usepackage{amsfonts}
\usepackage{amsmath}
\usepackage{amssymb}
\usepackage{eurosym}
\usepackage{mathrsfs}
\usepackage{palatino}
\usepackage{color}
\usepackage{esint}
\usepackage{url}
 \usepackage{color}

 \usepackage{amsrefs}     
 \usepackage[normalem]{ulem}
 \usepackage{hyperref} 
\hypersetup{
    colorlinks=false,       
    linkcolor=blue,          
    citecolor=magenta,        
    filecolor=magenta,      
    urlcolor=cyan           
}

\newcommand{\R}{\mathbb{R}}

\newcommand{\N}{\mathbb{N}}

\newcommand{\Z}{\mathbb{Z}}

\numberwithin{equation}{section}
\usepackage{graphicx}

\swapnumbers
\theoremstyle{plain}
\newtheorem{thm}[equation]{Theorem}
\newtheorem{lem}[equation]{Lemma}
\newtheorem{prop}[equation]{Proposition}

\theoremstyle{definition}
\newtheorem{defn}[equation]{Definition}

\theoremstyle{remark}
\newtheorem{rem}[equation]{Remark}

\author{Michael T. Lacey}
\address[M.T.L.]{ School of Mathematics, Georgia Institute of Technology, Atlanta GA 30332, USA}
\email {lacey@math.gatech.edu}
\thanks{Research of M.T.L. is supported in part by grant NSF-DMS 0968499, and the Australian Research Council through grant ARC-DP120100399.}

\author{Henri Martikainen}
\address[H.M.]{Department of Mathematics and Statistics, University of Helsinki, P.O.B. 68, FI-00014 Helsinki, Finland}
\email{henri.martikainen@helsinki.fi}
\thanks{This paper was completed while H.M. was still at Universit\'e Paris-Sud 11, Orsay. During this period the research of H.M. was supported by the Emil Aaltonen Foundation. }

\makeatletter
\@namedef{subjclassname@2010}{%
  \textup{2010} Mathematics Subject Classification}
\makeatother


\title[Local $Tb$ theorem with $L^2$ testing conditions and general measures]{Local $Tb$ theorem with $L^2$ testing conditions and general measures: Square functions}

\begin{document}
\begin{abstract}
Local $Tb$ theorems with $L^p$ type testing conditions, which are not scale invariant, have been studied widely in the case of the Lebesgue measure. In the non-homogeneous world
local $Tb$ theorems have only been proved assuming scale invariant ($L^{\infty}$ or BMO) testing conditions. In this paper, for the first time, we overcome these obstacles in the non-homogeneous world,
and prove a non-homogeneous local $Tb$ theorem with $L^2$ type testing conditions. This paper is in the setting of the vertical and conical square functions defined using general measures and kernels.
The proof uses various recent innovations including a Whitney averaging formula and the insertion of a Calder\'on--Zygmund stopping data of a fixed function
in to the construction of the twisted martingale difference operators.
\end{abstract}

\maketitle
\section{Introduction}
We assume that the Borel measure $\mu$ and the exponent $m > 0$ are related by the condition
\begin{equation*}\label{powerBound} \mu(B(x,r)) \lesssim r^{m}, \qquad x \in \R^{n}, \; r > 0. \end{equation*}
We are also given linear operators $\theta_{t}$, $t > 0$, which have the form
\begin{displaymath} \theta_{t}f(x) = \int_{\R^{n}} s_{t}(x,y)f(y) \, d\mu(y). \end{displaymath}
The kernels $s_{t}$ satisfy the size and continuity conditions
\begin{equation}\label{eq:size}
|s_t(x,y)| \lesssim \frac{t^{\alpha}}{(t+|x-y|)^{m+\alpha}}
\end{equation}
and
\begin{equation}\label{eq:hol}
|s_t(x,y) - s_t(x,z)| \lesssim \frac{|y-z|^{\alpha}}{(t+|x-y|)^{m+\alpha}} \end{equation}
whenever $|y-z| < t/2$. Here $\alpha > 0$ is a fixed constant. Notice that no regularity from $s_{t}$ is required in the first variable.

We study the conical and vertical square functions operators $S$ and $V$, defined by
\begin{displaymath} Sf(x) = \Big( \iint_{\Gamma(x)} |\theta_t f(y)|^2 \,\frac{d\mu(y) dt}{t^{m+1}}\Big)^{1/2} \quad \text{and} \quad Vf(x) = \left( \int_{0}^{\infty} |\theta_{t}f(y)|^{2} \, \frac{dt}{t} \right)^{1/2}. \end{displaymath}
In the definition of $S$ the cone $\Gamma(x) = \{(y,t) \in \R^{n+1}_+\colon \, |x-y| < t\}$, $x \in \R^n$, appears.
Let us note that since we are interested in the $L^2(\mu)$ boundedness of these operators, it is actually enough to study only vertical square functions.
This is because $\|Sf\|_{L^{2}(\mu)} = \|\tilde{V}f\|_{L^{2}(\mu)}$, where $\tilde{V}$ is the vertical square function with kernel
\begin{displaymath} \tilde{s}_{t}(x,y) = \left(\frac{\mu(B(x,t))}{t^{m}}\right)^{1/2}s_{t}(x,y). \end{displaymath}
Since the $x$-continuity of the kernels is not required, the kernel $\tilde{s}_{t}$ remains in our framework i.e. satisfies the assumptions \eqref{eq:size} and \eqref{eq:hol}.

We are interested in very general type $L^2(\mu)$ boundedness characterizations for $S$ and $V$.
With this we mean a local $Tb$ theorem with $L^2$ type testing conditions. Such theorems
have been widely studied in the homogeneous situation. However, in the non-homogeneous situation no such theorems appear in the previous literature (for Calder\'on--Zygmund operators nor square functions).
Indeed, in the setting of general measures all the previous literature assumes scale invariant $L^{\infty}(\mu)$ or BMO$(\mu)$ type testing conditions. Here, for the first time, we overcome these restrictions
and prove a non-homogeneous local $Tb$ theorem with $L^2$ type testing conditions:
\begin{thm}\label{thm:main} Assume that to every cube $Q \subset \R^n$ there is associated a function $b_Q$ which satisfies:
\begin{enumerate}
\item spt$\,b_Q \subset Q$;
\item $\langle b_Q \rangle_Q = 1$;
\item $\|b_Q\|_{L^2(\mu)}^2 \lesssim \mu(Q)$;
\item \begin{displaymath}
\iint_{\widehat Q} |\theta_t b_Q(x)|^2 d\mu(x) \frac{dt}{t} \lesssim \mu(Q), \qquad \widehat Q = Q \times (0, \ell(Q)).
\end{displaymath}
\end{enumerate}
Then $V$ is bounded on $L^{2}(\mu)$, that is we have the square function estimate
\begin{equation}\label{eq:SFbound}
\iint_{\R^{n+1}_+} |\theta_t f(x)|^2\,d\mu(x)\frac{dt}{t} \lesssim \|f\|_{L^2(\mu)}^2, \qquad f \in L^2(\mu).
\end{equation}
\end{thm}

After stating our main result, let us now discuss the history and references in more detail.
First local $Tb$ theorem, with $L^{\infty}$ control of the test functions and their images, is by Christ \cite{Ch}. This was proven for doubling measures.
Nazarov, Treil and Volberg \cite{NTVa} obtained a non-homogeneous version of this theorem.

The idea of using (in the homogeneous situation) just local $L^p$ type testing conditions was introduced by Auscher, Hofmann, Muscalu, Tao and Thiele \cite{AHMTT}.
However, their proof works only for the so-called perfect dyadic singular integral operators.
The assumptions are of the form $\int_Q |b^1_Q|^p \le |Q|$, $\int_Q |b^2_Q|^q \le |Q|$, $\int_Q |Tb^1_Q|^{q'} \le |Q|$ and $\int_Q |T^*b^2_Q|^{p'} \le |Q|$, where $s'$ denotes the dual exponent of $s$ and $1 < p, q \le \infty$.

It turned out to be difficult to extend these theorems to the general Calder\'on--Zygmund operators -- at least with the full range of exponents $p,q \in (1,\infty]$ (or even $p=q=2$ which is our main interest here).
Hofmann \cite{Ho1} was able to extend to general Calder\'on--Zygmund operators but at the price of needing a stronger
set of assumptions: $\int_Q |b^1_Q|^s \le |Q|$, $\int_Q |b^2_Q|^s \le |Q|$, $\int_Q |Tb^1_Q|^{2} \le |Q|$ and $\int_Q |T^*b^2_Q|^{2} \le |Q|$ for some $s > 2$. Auscher and Yang \cite{AY}
established the theorem for standard Calder\'on--Zygmund operators in the case $1/p + 1/q \le 1$ (and thus in the case $p=q=2$).

This chase for the most general range of exponents is not our main focus in this paper. We content on establishing the non-homogeneous result in the already very interesting case of $p=q=2$. But we still mention that a lot of further work
has been done to improve these exponents in the homogeneous situation. Hofmann \cite{Ho2} has given a full solution in the case of square functions. In the Calder\'on--Zygmund world the work of Auscher and Routin \cite{AR}
continued to shed some light to the general case of exponents, however, not giving a definite answer and involving additional technical conditions. The (almost) full solution is given by Hyt\"onen and Nazarov \cite{HN}.

We work with square functions (instead of Calder\'on--Zygmund operators) in this paper, since it allows a technically clearer framework. Indeed, it completely circumvents certain technicalities which arise with probabilistic methods in the
non-homogeneous local Tb situation in the case of Calder\'on--Zygmund operators. This has been elaborated in Remark 4.1 \cite{HM} and in \cite{LV-direct} (see especially Remark 2.14). In a following paper we intent to push these results  
to the Calder\'on--Zygmund case. But, like said, here we want to focus on the already critical problems that arise from using $L^2$ testing conditions with general measures.

Next, we discuss the proof. The proof is started by using the recent averaging identity over good Whitney regions by Martikainen and Mourgoglou \cite{MM}.
Such an identity is inspired by Hyt\"onen's proof of the $A_2$ conjecture \cite{Hy}, which uses a very nice refinement of the Nazarov--Treil--Volberg method of random dyadic systems.

The paper \cite{MM} gave a non-homogeneous global $Tb$ theorem for square functions. One result of the larger paper \cite{MMO} by Martikainen, Mourgoglou and Orponen is an extension of this result to the local situation. However, only with scale invariant testing
conditions. While using this proof as a general guide (we still give all the details)
we need new ideas to work with our much more general test functions. We may use some tricks from Hyt\"onen--Martikainen \cite{HM}. However, severe problems with general measures arise.
To circumvent these, the main idea is to insert a certain Calder\'on--Zygmund stopping data of a fixed function
in to the construction of the twisted martingale difference operators. This is inspired by the school of ideas by Lacey-Sawyer-Shen-Uriarte-Tuero \cite{LSUT} 
in the two weight Hilbert situation, and Lacey--V\"ah\"akangas \cites{LV-perfect,LV-dual}
in some related $Tb$ settings.  (See the comprehensive survey of the two weight Hilbert inequality \cite{primer}.) 

To conclude, we remark that the main theorem, Theorem \ref{thm:main}, can be proved assuming only that $\mu(B(x,r)) \le \lambda(x,r)$ for some $\lambda\colon \R^n \times (0,\infty) \to (0,\infty)$ satisfying that
$r \mapsto \lambda(x,r)$ is non-decreasing and $\lambda(x, 2r) \le C_{\lambda}\lambda(x,r)$ for all $x \in \R^n$ and $r > 0$. In this case one only needs to replace the kernel estimates by
\begin{displaymath}
|s_t(x,y)| \lesssim \frac{t^{\alpha}}{t^{\alpha}\lambda(x,t) + |x-y|^{\alpha}\lambda(x, |x-y|)}
\end{displaymath}
and
\begin{displaymath}
|s_t(x,y) - s_t(x,z)| \lesssim \frac{|y-z|^{\alpha}}{t^{\alpha}\lambda(x,t) + |x-y|^{\alpha}\lambda(x, |x-y|)}
\end{displaymath}
whenever $|y-z| < t/2$. This is done in the global situation in \cite{MM}. Here we skip the required modifications.

\subsection{Notation}
We consider a dyadic grid $\mathcal{D}$ in $\R^n$. For $Q \in \mathcal{D}$ we have:
\begin{itemize}
\item $\ell(Q)$ is the side length of $Q$;
\item $\widehat Q = Q \times (0, \ell(Q))$ is the Carleson box associated with $Q$;
\item $W_Q = Q \times (\ell(Q)/2, \ell(Q))$ is the Whitney region associated with $Q$;
\item ch$(Q) = \{Q' \in \mathcal{D}:\, Q' \subset Q, \, \ell(Q') = \ell(Q)/2\}$;
\item gen$(Q)$ is determined by $\ell(Q) = 2^{-\textup{gen}(Q)}$;
\item $Q^{(k)} \in \mathcal{D}$ is the unique cube for which $\ell(Q^{(k)}) = 2^k\ell(Q)$ and $Q \subset Q^{(k)}$.
\end{itemize}

\section{Twisted martingale difference operators $\Delta_Q$}
Consider a fixed dyadic grid $\mathcal{D}$ in $\R^n$.
Let $f$ be a fixed function and let $Q^* \in \mathcal{D}$ be a fixed dyadic cube. Set $\mathcal{F}^0_{Q^*} = \{Q^*\}$ and $\alpha_f(Q^*) = \langle |f| \rangle_{Q^*}$. Let $\mathcal{F}_{Q^*}^1$ consist of the maximal cubes
$Q \in \mathcal{D}$, $Q \subset Q^*$, for which at least one of the following three conditions holds:
\begin{enumerate}
\item $|\langle b_{Q^*} \rangle_Q| < 1/2$;
\item $\langle |b_{Q*}|^2 \rangle_Q > 16A^2$;
\item $\langle |f| \rangle_Q > 32A \cdot \alpha_f(Q^*)$. 
\end{enumerate}
Here $A$ is the implied constant of the assumption (3) of Theorem \ref{thm:main}: $\|b_R\|_{L^2(\mu)}^2 \le A\mu(R)$ for every cube $R \subset \R^n$.
We define
\begin{displaymath}
\alpha_f(Q) = \left\{ \begin{array}{ll}
\alpha_f(Q^*) & \textrm{if } \langle |f|\rangle_Q < 2\alpha_f(Q^*), \\
\langle |f|\rangle_Q  & \textrm{if }  \langle |f|\rangle_Q  \ge 2\alpha_f(Q^*). \end{array} \right.
\end{displaymath}
Interestingly, the stopping time (3) is so potent that we do not even need to perform a stopping with respect to the condition (4) of Theorem \ref{thm:main}!

Next, one repeats the previous procedure by replacing $Q^*$ with a fixed $Q \in \mathcal{F}^1_{Q^*}$. The combined collection of stopping
cubes resulting from this is called $\mathcal{F}^2_{Q^*}$. This is continued and one sets $\mathcal{F}_{Q^*} = \bigcup_{j=0}^{\infty} \mathcal{F}^j_{Q^*}$.
Finally, for every $Q \in \mathcal{D}$, $Q \subset Q^*$, we let $Q^a \in \mathcal{F}_{Q^*}$ be the minimal
cube $R \in \mathcal{F}_{Q^*}$ for which $Q \subset R$.

\begin{lem}
We have the following Carleson estimate:
\begin{displaymath}
\mathop{\sum_{F \in \mathcal{F}_{Q^*}}}_{F \subset Q} \mu(F) \lesssim \mu(Q), \qquad Q \in \mathcal{D}, \, Q \subset Q^*.
\end{displaymath}
\end{lem}
\begin{proof}
Let $F \in \mathcal{F}_{Q^*}$. Consider a disjoint collection of dyadic cubes $\{Q^1_i\}_i$ for which $Q^1_i \subset F$ and $|\langle b_F \rangle_{Q^1_i}| < 1/2$.
We have that
\begin{align*}
\mu(F) = \int_F b_F \,d\mu &= \int_{F \setminus \bigcup_i Q^1_i} b_F\,d\mu + \sum_i \int_{Q^1_i} b_F \,d\mu \\
&\le \mu\Big(F \setminus \bigcup_i Q^1_i\Big)^{1/2} \Big( \int_{F} |b_F|^2 \,d\mu \Big)^{1/2} + \frac{1}{2} \sum_i \mu(Q^1_i) \\
&\le A^{1/2} \mu\Big(F \setminus \bigcup_i Q^1_i\Big)^{1/2} \mu(F)^{1/2} + \frac{1}{2} \mu(F),
\end{align*}
which implies that
\begin{displaymath}
\mu(F) \le 4A \cdot \mu\Big(F \setminus \bigcup_i Q^1_i\Big) = 4A\Big[ \mu(F) - \mu\Big(\bigcup_i Q^1_i\Big)\Big].
\end{displaymath}
From here we can read that
\begin{displaymath}
\mu\Big(\bigcup_i Q^1_i\Big) \le \Big(1- \frac{1}{4A}\Big)\mu(F).
\end{displaymath}

Consider then a disjoint collection of dyadic cubes $\{Q^2_i\}_i$ for which $Q^2_i \subset F$ and $\langle |b_F|^2 \rangle_{Q^2_i} > 16A^2$. Simply notice that now
\begin{displaymath}
\mu\Big(\bigcup_i Q^2_i\Big) \le \frac{1}{16A^2} \sum_i \int_{Q^2_i} |b_F|^2\,d\mu \le \frac{1}{16A^2} \int_{F} |b_F|^2\,d\mu \le \frac{1}{16A} \mu(F).
\end{displaymath}

Lastly, consider a disjoint collection of dyadic cubes $\{Q^3_i\}_i$ for which $Q^3_i \subset F$ and $\langle |f| \rangle_{Q^3_i} > 32A \cdot \alpha_f(F)$.
If $F \in \mathcal{F}^j$, let $S \in \mathcal{F}^{j-1}$ be such that $F \subset S$. We can have $\alpha_f(F) = \langle |f|\rangle_F$, or
$\alpha_f(F) = \alpha_f(S)$. But in the latter case we must have $\langle |f| \rangle_F < 2\alpha_f(S) = 2\alpha_f(F)$. Therefore, we always have
$\langle |f| \rangle_{Q^3_i} \ge 16A \langle |f|\rangle_F$. From this it is immediately clear that
\begin{displaymath}
\mu\Big(\bigcup_i Q^3_i\Big) \le \frac{1}{16A} \mu(F).
\end{displaymath}

Combining the analysis we may conclude that for $F \in \mathcal{F}^j_{Q^*}$ there holds that
\begin{displaymath}
\mathop{\sum_{S \in \mathcal{F}^{j+1}_{Q^*}}}_{S \subset F} \mu(S) \le \Big(1 - \frac{1}{8A}\Big)\mu(F) =: \tau\mu(F), \qquad \tau < 1.
\end{displaymath}
From this it is easy to conclude that if $Q \in \mathcal{D}$, $Q \subset Q^*$ and $Q^a \in \mathcal{F}^j_{Q^*}$, then there holds that
\begin{displaymath}
\mathop{\sum_{F \in \mathcal{F}^{j+k}_{Q^*}}}_{F \subset Q} \mu(F) \le \tau^{k-1}\mu(Q), \qquad k \ge 1.
\end{displaymath}
Finally, this implies that
\begin{displaymath}
\mathop{\sum_{F \in \mathcal{F}_{Q^*}}}_{F \subset Q} \mu(F) \le \Big(1 + \frac{1}{1-\tau}\Big)\mu(Q) = (1+8A)\mu(Q).
\end{displaymath}

\end{proof}
\begin{lem}
We have the Calder\'on--Zygmund stopping data estimate for the fixed function $f$:
\begin{displaymath}
\sum_{F \in \mathcal{F}_{Q^*}} \alpha_f(F)^2\mu(F) \lesssim \|f\|_{L^2(Q^*; \mu)}^2.
\end{displaymath}
\end{lem}
\begin{proof}
We define
\begin{displaymath}
\mathcal{G} := \{Q^*\} \cup \bigcup_{j=0}^{\infty} \bigcup_{S \in \mathcal{F}^j_{Q^*}} \{F \in \mathcal{F}^{j+1}_{Q^*}: \, F \subset S \textrm{ and } \langle |f| \rangle_F \ge 2\alpha_f(S)\}.
\end{displaymath}
For every $F \in \mathcal{F}$ we let $F^b \in \mathcal{G}$ be the minimal cube $R \in \mathcal{G}$ for which $F \subset R$. We have that
\begin{displaymath}
\sum_{F \in \mathcal{F}_{Q^*}} \alpha_f(F)^2\mu(F) = \sum_{G \in \mathcal{G}} \alpha_f(G)^2 \mathop{\sum_{F \in \mathcal{F}}}_{F^b = G} \mu(F) \lesssim \sum_{G \in \mathcal{G}} \alpha_f(G)^2 \mu(G).
\end{displaymath}
Note that $\mu$-a.e. $x \in Q^*$ belongs to only finitely many $G \in \mathcal{G}$. Let $x$ be such and let $S = S_x \in \mathcal{G}$ be the minimal $R \in \mathcal{G}$ for which $x \in R$. We have that
\begin{displaymath}
\sum_{G \in \mathcal{G}} \alpha_f(G)^2 1_G(x) = \sum_{G \in \mathcal{G}} \langle |f| \rangle_G^2 1_G(x) \le \Big( \sum_{j=0}^{\infty} 2^{-j}\Big) \langle |f| \rangle_S^2 \lesssim M_{\mathcal{D}}^{\mu} (f1_{Q^*})(x)^2,
\end{displaymath}
where $M_{\mathcal{D}}^{\mu}$ is the dyadic maximal function $M_{\mathcal{D}}^{\mu} g(x) := \sup_{Q \in \mathcal{D}} 1_Q(x) \langle |g| \rangle_Q$. The claim readily follows.
\end{proof}
\begin{defn}
If $Q \in \mathcal{D}$, $Q \subset Q^*$, we define the twisted martingale difference operators
\begin{equation*}
\Delta_Q g = \sum_{Q' \in \, \textrm{ch}(Q)} \Big[\frac{\langle g \rangle_{Q'}}{\langle b_{(Q')^a}\rangle_{Q'}}b_{(Q')^a} - \frac{\langle g \rangle_Q}{\langle b_{Q^a}\rangle_Q}b_{Q^a}\Big]1_{Q'}.
\end{equation*}
\end{defn}
\begin{rem}
Notice that we have included the Calder\'on--Zygmund stopping data (4) of the fixed function $f$. This means that the twisted martingale difference operators $\Delta_Q = \Delta_Q^f$, unlike usually, depend also on the fixed function $f$.
This is key to be able to prove the following square function estimate in the non-homogeneous situation!
\end{rem}

\begin{prop}\label{prop:SF}
For the fixed function $f$ we have the estimate
\begin{displaymath}
\mathop{\sum_{Q \in \mathcal{D}}}_{Q \subset Q^*} \|\Delta_Q f\|_{L^2(\mu)}^2 \lesssim \|f\|_{L^2(Q^*; \mu)}^2.
\end{displaymath}
\end{prop}
\begin{proof}
Let us write
\begin{displaymath}
\mathop{\sum_{Q \in \mathcal{D}}}_{Q \subset Q^*} \|\Delta_Q f\|_{L^2(\mu)}^2 = A + B,
\end{displaymath}
where
\begin{align*}
A &:= \mathop{\sum_{Q \in \mathcal{D}}}_{Q \subset Q^*} \mathop{\sum_{Q' \in \, \textrm{ch}(Q)}}_{(Q')^a = Q'} \int_{Q'} \Big| 
\frac{\langle f \rangle_{Q'}}{\langle b_{Q'}\rangle_{Q'}}b_{Q'} - \frac{\langle f \rangle_Q}{\langle b_{Q^a}\rangle_Q}b_{Q^a}\Big|^2\,d\mu; \\
B &:= \mathop{\sum_{Q \in \mathcal{D}}}_{Q \subset Q^*} \mathop{\sum_{Q' \in \, \textrm{ch}(Q)}}_{(Q')^a = Q^a}  \Big|  \frac{\langle f \rangle_{Q'}}{\langle b_{Q^a}\rangle_{Q'}} - \frac{\langle f \rangle_Q}{\langle b_{Q^a}\rangle_Q} \Big|^2
\int_{Q'} |b_{(Q')^a}|^2\,d\mu.
\end{align*}

Notice first that
\begin{equation}\label{Bterm}
B \lesssim \mathop{\sum_{Q \in \mathcal{D}}}_{Q \subset Q^*} \mathop{\sum_{Q' \in \, \textrm{ch}(Q)}}_{(Q')^a = Q^a}[ |\langle f \rangle_{Q'}|^2 |\langle b_{Q^a} \rangle_{Q'} - \langle b_{Q^a} \rangle_{Q}|^2 \mu(Q')
+ |\langle f \rangle_{Q'} - \langle f \rangle_Q|^2 \mu(Q')].
\end{equation}
The first term in \eqref{Bterm} is controlled by checking that the sequence
\begin{displaymath}
|\langle b_{Q^a} \rangle_Q - \langle b_{Q^a} \rangle_{Q^{(1)}}|^2 \mu(Q), \qquad Q \in \mathcal{D}, \, Q \subset Q^*,
\end{displaymath}
is Carleson. For this, just note that:
\begin{align*}
&K(Q) := \sum_{S \subset Q}  |\langle b_{S^a} \rangle_S - \langle b_{S^a} \rangle_{S^{(1)}}|^2\mu(S) =  |\langle b_{Q^a} \rangle_Q - \langle b_{Q^a} \rangle_{Q^{(1)}}|^2\mu(Q) \\
&+ \mathop{\sum_{S \subsetneq Q}}_{S^a = Q^a} |\langle b_{Q^a}1_Q \rangle_S - \langle b_{Q^a}1_Q \rangle_{S^{(1)}}|^2\mu(S)  
+ \mathop{\sum_{H \subsetneq Q}}_{H^a = H} \sum_{S:\,S^a = H} |\langle b_{H} \rangle_S - \langle b_{H} \rangle_{S^{(1)}}|^2\mu(S).
\end{align*}
In the very first term (with $S =Q $) notice that indeed $|\langle b_{Q^a} \rangle_{Q^{(1)}}|^2 \lesssim 1$ (if $Q = Q^a$, then $\int_{Q^{(1)}} |b_{Q^a}|^2\,d\mu = \int_{Q} |b_{Q}|^2\,d\mu \lesssim \mu(Q) \le \mu(Q^{(1)})$, and otherwise
clearly $\int_{Q^{(1)}} |b_{Q^a}|^2\,d\mu \lesssim \mu(Q^{(1)}))$. Using the unweighted square function estimate we thus get that
\begin{displaymath}
K(Q) \lesssim \mu(Q) + \|b_{Q^a}1_Q\|_{L^2(\mu)}^2 +  \mathop{\sum_{H \subset Q}}_{H^a = H} \|b_H\|_{L^2(\mu)}^2 \lesssim \mu(Q) +  \mathop{\sum_{H \subset Q}}_{H^a = H} \mu(H) \lesssim \mu(Q).
\end{displaymath}

The second term in \eqref{Bterm} is simply controlled by the unweighted square function estimate. Therefore, we get that $B \lesssim \|f\|_{L^2(Q^*; \mu)}^2$.

Next, notice that there holds that
\begin{displaymath}
\mathop{\sum_{Q \in \mathcal{D}}}_{Q \subset Q^*} \mathop{\sum_{Q' \in \, \textrm{ch}(Q)}}_{(Q')^a = Q'} \frac{|\langle f \rangle_{Q'}|^2}{|\langle b_{Q'}\rangle_{Q'}|^2}\int_{Q'} |b_{Q'}|^2\,d\mu \lesssim 
\mathop{\sum_{R:\, R^a = R}}_{R \subset Q^*} |\langle f \rangle_R|^2 \mu(R) \lesssim  \|f\|_{L^2(Q^*; \mu)}^2.
\end{displaymath}
These estimates have not yet used the fact that we included the Calder\'on--Zygmund stopping data (4) of $f$ i.e. they would hold for every function.

To control $A$, it only remains to control
\begin{displaymath}
C := \mathop{\sum_{Q \in \mathcal{D}}}_{Q \subset Q^*} \mathop{\sum_{Q' \in \, \textrm{ch}(Q)}}_{(Q')^a = Q'} |\langle f \rangle_Q|^2 \int_{Q'} |b_{Q^a}|^2\,d\mu.
\end{displaymath}
For this we write
\begin{displaymath}
C = \sum_{F \in \mathcal{F}_{Q^*}} \sum_{Q: \, Q^a = F}  \mathop{\sum_{Q' \in \, \textrm{ch}(Q)}}_{(Q')^a = Q'} |\langle f \rangle_Q|^2 \int_{Q'} |b_{F}|^2\,d\mu =: \sum_{F \in \mathcal{F}} C_F.
\end{displaymath}
Utilizing the Calder\'on--Zygmund stopping data (4) of $f$ here, it is enough to show that $C_F \lesssim \alpha_f(F)^2\mu(F)$. But we have that
\begin{displaymath}
C_F \lesssim \alpha_f(F)^2  \sum_{Q: \, Q^a = F}  \mathop{\sum_{Q' \in \, \textrm{ch}(Q)}}_{(Q')^a = Q'}  \int_{Q'} |b_{F}|^2\,d\mu \le \alpha_f(F)^2   \int_{F} |b_{F}|^2\,d\mu \lesssim \alpha_f(F)^2\mu(F).
\end{displaymath}
Here we used that the stopping cubes $Q' \subset F$ are disjoint since they are of the same generation (because of the condition $Q^a = F$).
\end{proof}

\begin{rem}
At least in the doubling situation (the measure $\mu$ would be doubling) this estimate can be made uniform (to work for every function). Indeed, then one can bound
\begin{displaymath}
\int_{Q'} |b_{Q^a}|^2\,d\mu \le \int_{Q} |b_{Q^a}|^2\,d\mu \lesssim \mu(Q) \lesssim \mu(Q'),
\end{displaymath}
after which one can bound
\begin{displaymath}
C \lesssim  \mathop{\sum_{Q \in \mathcal{D}}}_{Q \subset Q^*} \mathop{\sum_{Q' \in \, \textrm{ch}(Q)}}_{(Q')^a = Q'}[ |\langle f \rangle_{Q'} - \langle f \rangle_Q|^2 \mu(Q')
+ |\langle f \rangle_{Q'}|^2 \mu(Q')].
\end{displaymath}
Here the first term can be bounded by the unweighted square function estimate, and the second by using the Carleson property of the stopping cubes. Therefore, in the doubling
situation all the needed estimates can be completed without using the stopping data of $f$.
\end{rem}

\begin{prop}
One can write $\mu$-a.e. and in $L^2(\mu)$ that
\begin{displaymath}
g1_{Q^*} = \mathop{\sum_{Q \in \mathcal{D}}}_{Q \subset Q^*} \Delta_Q g + \langle g \rangle_{Q^*} b_{Q^*}
\end{displaymath}
for any function $g \in L^2(\mu)$.
\end{prop}
\begin{proof}
The pointwise identity is a standard deduction. The $L^2(\mu)$ identity then follows by dominated convergence using the fact
that
\begin{displaymath}
\Big( \sum_{Q:\, Q^a = Q} |b_Q|^2 \Big)^{1/2} \in L^2(\mu).
\end{displaymath}
\end{proof}

\section{Beginning of the proof of the local $Tb$}
\subsection{Random dyadic grids}
We need to insert a standard disclaimer about random dyadic grids.
To this end, let us be given a random dyadic grid $\mathcal{D} = \mathcal{D}(w)$, $w = (w_i)_{i \in \Z} \in (\{0,1\}^n)^{\Z}$.
This means that $\mathcal{D} = \{Q + \sum_{i:\, 2^{-i} < \ell(Q)} 2^{-i}w_i: \, Q \in \mathcal{D}_0\} = \{Q + w: \, Q \in \mathcal{D}_0\}$, where we simply have defined
$Q + w := Q + \sum_{i:\, 2^{-i} < \ell(Q)} 2^{-i}w_i$. Here $\mathcal{D}_0$ is the standard dyadic grid of $\R^n$.

We set $\gamma = \alpha/(2m+2\alpha)$, where $\alpha > 0$ appears in the kernel estimates and $m$ appears in $\mu(B(x,r)) \lesssim r^m$.
A cube $Q \in \mathcal{D}$ is called bad if there exists another cube $\tilde Q \in \mathcal{D}$ so that $\ell(\tilde Q) \ge 2^r \ell(Q)$ and $d(Q, \partial \tilde Q) \le \ell(Q)^{\gamma}\ell(\tilde Q)^{1-\gamma}$.
Otherwise it is good.
One notes that $\pi_{\textrm{good}} := \mathbb{P}_{w}(Q + w \textrm{ is good})$ is independent of $Q \in \mathcal{D}_0$. The parameter $r$ is a fixed constant so large that $\pi_{\textrm{good}} > 0$
and $2^{r(1-\gamma)} \ge 3$.

Furthermore, it is important to note that for a fixed $Q \in \mathcal{D}_0$
the set $Q + w$ depends on $w_i$ with $2^{-i} < \ell(Q)$, while the goodness (or badness) of $Q + w$ depends on $w_i$ with $2^{-i} \ge \ell(Q)$. In particular, these notions are independent (meaning that
for any fixed $Q \in \mathcal{D}_0$ the random variable $w \mapsto 1_{\textup{good}}(Q+w)$ and any random variable that depends only on the cube $Q+w$ as a set, like $w \mapsto \int_{Q+w} f\,d\mu$, are independent).

\subsection{Whitney averaging identity}
Fix a compactly supported function $f$. We will prove \eqref{eq:SFbound} for this function.
The proof is started by writing the identity
\begin{align*}
\iint_{\R^{n+1}_+} |\theta_t f(x)|^2\,d\mu(x)\frac{dt}{t} = 
\frac{1}{\pi_{\textrm{good}}} E_w \sum_{R \in \mathcal{D}(w)_{\textup{good}}} \iint_{W_R} |\theta_t f(x)|^2\,d\mu(x)\frac{dt}{t},
\end{align*}
where $W_R = R \times (\ell(R)/2, \ell(R))$ is the Whitney region associated with $R \in \mathcal{D} = \mathcal{D}(w)$. The proof of the identity is based on the facts that for every fixed $R \in \mathcal{D}_0$
the random variables $1_{\textup{good}}(R+w)$ and $\iint_{W_{R+w}} |\theta_t f(x)|^2\,d\mu(x)\frac{dt}{t}$ are independent, and that we have $\pi_{\textrm{good}} = \mathbb{P}_{w}(R + w \textrm{ is good}) = E_w 1_{\textup{good}}(R+w)$.

We fix the grid $\mathcal{D} = \mathcal{D}(w)$ i.e. we fix $w$ from the probability space.
It is enough to prove that for any fixed large $s$ there holds that
\begin{displaymath}
\mathop{\sum_{R \in \mathcal{D}_{\textup{good}}}}_{\ell(R) \le 2^s}  \iint_{W_R} |\theta_t f(x)|^2\,d\mu(x)\frac{dt}{t} \lesssim \|f\|_{L^2(\mu)}^2.
\end{displaymath}

Now fix $N \in \N$ such that spt$\ f \subset B(0,2^N)$ and consider any $s \ge N$. We shall expand $f$ using the twisted martingale difference operators $\Delta_Q = \Delta_Q^f$ with stopping data adapted to $f$.
Write
\begin{displaymath}
f = \mathop{\mathop{\sum_{Q^* \in \mathcal{D}}}_{\ell(Q^*) = 2^s}}_{Q^* \cap B(0, 2^N) \ne \emptyset} \mathop{\sum_{Q \in \mathcal{D}}}_{Q \subset Q^*} \Delta_Q f.
\end{displaymath}
On the largest $Q^*$ level we agree (by abuse of notation) that $\Delta_{Q^*} = E_{Q^*} + \Delta_{Q^*}$, where
$E_{Q^*}f = \langle f \rangle_{Q^*}b_{Q^*}$. Therefore, we have that $\int \Delta_Qf \,d\mu = 0$ except when $Q  = Q^*$ for some $Q^*$ with $\ell(Q^*) = 2^s$.
Since $\#\{Q^* \in \mathcal{D}:\, Q^* \cap B(0,2^N) \ne \emptyset\} \lesssim 1$, we can fix one $Q^*$ with $\ell(Q^*) = 2^s$, and concentrate on proving that
\begin{displaymath}
\mathop{\sum_{R \in \mathcal{D}_{\textup{good}}}}_{\ell(R) \le 2^s}  \iint_{W_R} \Big|\mathop{\sum_{Q \in \mathcal{D}}}_{Q \subset Q^*}\theta_t \Delta_Q f(x)\Big|^2\,d\mu(x)\frac{dt}{t} \lesssim \|f\|_{L^2(\mu)}^2.
\end{displaymath}
This is done case by case in the chapters that follow. Many of the considerations below are of a standard nature. The special form of the stopping data eases the control of the paraproduct term.
The validity of the square function estimate for the function $f$, Proposition \ref{prop:SF}, is of key importance.

\section{The case $\ell(Q) < \ell(R)$}
Since $\ell(Q) < \ell(R) \le 2^s$, we have that $\int \Delta_Q f\,d\mu = 0$. Using this we write
\begin{align*}
 \theta_t \Delta_Q f(x) = \int_Q [s_t(x,y) - s_t(x, c_Q)]\Delta_Qf(y)\,d\mu(y), \qquad (x,t) \in W_R.
\end{align*}
Because $|y-c_Q| \le \ell(Q)/2 \le \ell(R)/4 < t/2$ for every $y \in Q$, we get using the full power of \eqref{eq:hol} that
\begin{displaymath}
|\theta_t \Delta_Q f(x)| \lesssim A_{QR}\mu(R)^{-1/2}\|\Delta_Qf\|_{L^2(\mu)}, \qquad (x,t) \in W_R.
\end{displaymath}
Here
\begin{align*}
A_{QR} &:= \frac{\ell(Q)^{\alpha/2}\ell(R)^{\alpha/2}}{D(Q,R)^{m+\alpha}} \mu(Q)^{1/2}\mu(R)^{1/2}; \\
D(Q,R) &:= \ell(Q) + \ell(R) + d(Q,R).
\end{align*}
Moreover, this $ \ell ^2 $ estimate holds 
\begin{displaymath}
\sum_R \Big[ \sum_Q A_{QR} x_Q \Big]^2 \lesssim \sum_Q x_Q^2.
\end{displaymath}
This yields that
\begin{align*}
\sum_{R:\, \ell(R) \le 2^s} \iint_{W_R} \Big| \mathop{\sum_{Q:\,Q \subset Q^*}}_{\ell(Q) < \ell(R)} \theta_t \Delta_Q f(x)\Big|^2\,d\mu(x)\frac{dt}{t} &\lesssim \sum_R \Big[  \sum_Q A_{QR} \|\Delta_Qf\|_{L^2(\mu)} \Big]^2 \\
&\lesssim \sum_Q \|\Delta_Q f\|_{L^2(\mu)}^2 \lesssim \|f\|_{L^2(\mu)}^2.
\end{align*}
\section{The case $\ell(Q) \ge \ell(R)$ and $d(Q,R) > \ell(R)^{\gamma}\ell(Q)^{1-\gamma}$}
The size estimate \eqref{eq:size} gives that
\begin{displaymath}
|\theta_t \Delta_Q f(x)| \lesssim  \frac{\ell(R)^{\alpha}}{d(Q,R)^{m+\alpha}}\mu(Q)^{1/2}\|\Delta_Qf\|_{L^2(\mu)}, \qquad (x,t) \in W_R.
\end{displaymath}
By the previous section it is enough to note that
\begin{equation}\label{eq:sep}
 \frac{\ell(R)^{\alpha}}{d(Q,R)^{m+\alpha}}\mu(Q)^{1/2} \lesssim A_{QR}\mu(R)^{-1/2}.
\end{equation}
In the case $d(Q,R) \ge \ell(Q)$ this is trivial. In the opposite case we note that $d(Q,R)^{m+\alpha} \gtrsim D(Q,R)^{m+\alpha}\ell(Q)^{-\alpha/2}\ell(R)^{\alpha/2}$.
This is seen by combining the facts that $d(Q,R) > \ell(R)^{\gamma}\ell(Q)^{1-\gamma}$, $\gamma m + \gamma \alpha = \alpha/2$ and $D(Q,R) \lesssim \ell(Q)$. This proves \eqref{eq:sep}
and therefore completes this section.

\section{The case $\ell(R) \le \ell(Q) \le 2^r\ell(R)$ and $d(Q,R) \le \ell(R)^{\gamma}\ell(Q)^{1-\gamma}$}
Begin by noting that
\begin{align*}
\sum_{R:\, \ell(R) \le 2^s} & \iint_{W_R} \Big| \mathop{\sum_{Q \subset Q^*:\, \ell(R) \le \ell(Q) \le  2^r\ell(R)}}_{d(Q,R) \le \ell(R)^{\gamma}\ell(Q)^{1-\gamma}} \theta_t \Delta_Q f(x)\Big|^2\,d\mu(x)\frac{dt}{t} \\
&\lesssim  \sum_{Q \subset Q^*} \sum_{R:\, R \sim Q}  \iint_{W_R} |\theta_t \Delta_Q f(x)|^2\,d\mu(x)\frac{dt}{t},
\end{align*}
where we have written $Q \sim R$ to mean $\ell(Q) \sim \ell(R)$ and $d(Q,R) \lesssim \min(\ell(Q), \ell(R))$.
We also used the fact that given $R$ there are $\lesssim 1$ cubes $Q$ for which $Q \sim R$.

Recalling that $\ell(Q) \sim \ell(R)$, the size estimate \eqref{eq:size} gives that
\begin{displaymath}
|\theta_t \Delta_Q f(x)| \lesssim t^{-m} \mu(Q)^{1/2} \|\Delta_Qf\|_{L^2(\mu)} \lesssim \mu(R)^{-1/2} \|\Delta_Qf\|_{L^2(\mu)}, \qquad (x,t) \in W_R.
\end{displaymath}
To complete this section, we note that this immediately gives that
\begin{displaymath}
\sum_{Q:\, Q \subset Q^*} \sum_{R:\, R \sim Q} \iint_{W_R} |\theta_t \Delta_Q f(x)|^2\,d\mu(x)\frac{dt}{t} \lesssim \sum_{Q:\, Q \subset Q^*}  \|\Delta_Q f\|_{L^2(\mu)}^2  \sum_{R:\, R \sim Q} 1 \lesssim \|f\|_{L^2(\mu)}^2.
\end{displaymath}

\section{The case $\ell(Q) > 2^r\ell(R)$ and $d(Q,R) \le \ell(R)^{\gamma}\ell(Q)^{1-\gamma}$}
The goodness of $R$ implies that here $R \subset Q$.
Therefore, we need to only consider the term
\begin{align*}\label{eq:maineq}
 \mathop{\mathop{\sum_{R \in \mathcal{D}_{\textup{good}}}}_{\ell(R) < 2^{s-r}}}_{R \subset Q^*} \iint_{W_R} \Big| \sum_{k=r+1}^{s+\textup{gen}(R)} \theta_t \Delta_{R^{(k)}} f(x)\Big|^2\,d\mu(x)\frac{dt}{t}.
\end{align*}
In what follows we shall not always write that everything is inside $Q^*$. We also write $\mathcal{F} = \mathcal{F}_{Q^*}$.

Before having to split the argument into a case study, we prove two lemmata which are useful in both cases.
\begin{lem}\label{lem:es1}
For $R \in \mathcal{D}_{\textup{good}}$ and $k \ge r+1$ there holds that
\begin{displaymath}
|\theta_t(1_{R^{(k)} \setminus R^{(k-1)}} \Delta_{R^{(k)}}f)(x)| \lesssim 2^{-\alpha k/2} \mu(R^{(k-1)})^{-1/2} \|\Delta_{R^{(k)}}f\|_{L^2(\mu)}, \qquad (x,t) \in W_R.
\end{displaymath}
\end{lem}
\begin{proof}
Let $S \in \textup{ch}(R^{(k)})$, $S \subset R^{(k)} \setminus R^{(k-1)}$.
The size estimate \eqref{eq:size} gives that
\begin{displaymath}
|\theta_t(1_S \Delta_{R^{(k)}}f)(x)| \lesssim \frac{\ell(R)^{\alpha}}{d(S,R)^{m+\alpha}} \int_S |\Delta_{R^{(k)}}f(y)|\,d\mu(y).
\end{displaymath}
The goodness gives that $d(S,R) > \ell(R)^{\gamma}\ell(S)^{1-\gamma}$, and thus
\begin{displaymath}
d(S,R)^{m+\alpha} > \ell(R)^{\alpha/2}\ell(S)^{\alpha/2} \ell(S)^m \gtrsim  \ell(R)^{\alpha/2}\ell(S)^{\alpha/2} \mu(S)^{1/2} \mu(R^{(k-1)})^{1/2}.
\end{displaymath}
A combination of these observations immediately yields the desired result.
\end{proof}
\begin{lem}\label{lem:es2}
For $R \in \mathcal{D}_{\textup{good}}$ and $k \ge r+1$ there holds that
\begin{displaymath}
|\theta_t(1_{(R^{(k-1)})^c} b_{(R^{(k)})^a}(x)| \lesssim 2^{-\alpha k/2}, \qquad (x,t) \in W_R.
\end{displaymath}
\end{lem}
\begin{proof}
Let $M \ge 0$ be such that $R^{(k+M)} = (R^{(k)})^a$. For $(x,t) \in W_R$ the size estimate \eqref{eq:size} gives that
\begin{align*}
|\theta_t(1_{(R^{(k-1)})^c} b_{(R^{(k)})^a}(x)| &\lesssim \int_{(R^{(k)})^a \setminus R^{(k-1)}} \frac{\ell(R)^{\alpha}}{|x-y|^{m+\alpha}} |b_{(R^{(k)})^a}(y)|\,d\mu(y) \\
&= \ell(R)^{\alpha} \sum_{j=0}^M \int_{R^{(k+j)} \setminus R^{(k+j-1)}} \frac{|b_{(R^{(k)})^a}(y)|}{|x-y|^{m+\alpha}}|\,d\mu(y).
\end{align*}
Notice that goodness gives that
\begin{align*}
|x-y|^{m+\alpha} &\ge d(R, \partial R^{(k+j-1)})^{m+\alpha} \\
&\ge \ell(R)^{\alpha/2} \ell(R^{(k+j-1)})^{m+\alpha/2} \gtrsim 2^{\alpha k/2} \ell(R)^{\alpha} 2^{\alpha j/2} \mu(R^{(k+j)}).
\end{align*}
Therefore, we have that
\begin{align*}
|\theta_t(1_{(R^{(k-1)})^c} b_{(R^{(k)})^a}(x)| \lesssim 2^{-\alpha k/2}  \sum_{j=0}^M 2^{-\alpha j /2}\mu(R^{(k+j)})^{-1}\int_{R^{(k+j)}} |b_{(R^{(k)})^a}| \,d\mu.
\end{align*}
This gives the desired result after noting that $\int_{R^{(k+j)}} |b_{(R^{(k)})^a}| \,d\mu \lesssim \mu(R^{(k+j)})$ for every $j = 0, 1, \ldots, M$.
\end{proof}

\subsection{The case $(R^{(k-1)})^a = (R^{(k)})^a$}
In this case we may write
\begin{align}\label{eq:split1}
\Delta_{R^{(k)}} f = 1_{R^{(k)} \setminus R^{(k-1)}} \Delta_{R^{(k)}}f - 1_{(R^{(k-1)})^c} B_{R^{(k-1)}} b_{(R^{(k)})^a} + B_{R^{(k-1)}}b_{(R^{(k)})^a},
\end{align}
where
\begin{displaymath}
B_{R^{(k-1)}} = \frac{\langle f \rangle_{R^{(k-1)}}}{\langle b_{(R^{(k-1)})^a}\rangle_{R^{(k-1)}}} - \frac{\langle f \rangle_{R^{(k)}} }{\langle b_{(R^{(k)})^a}\rangle_{R^{(k)} }}
\end{displaymath}
with the minus term missing if $\ell(Q) = 2^s$. Notice that the accretivity condition gives that
\begin{align*}
|B_{R^{(k-1)}}| \mu(R^{(k-1)}) \lesssim \Big| \int_{R^{(k-1)}} B_{R^{(k-1)}} b_{(R^{(k)})^a}\,d\mu\Big| &= \Big| \int_{R^{(k-1)}} \Delta_{R^{(k)}} f\,d\mu\Big|  \\
&\lesssim \mu(R^{(k-1)})^{1/2}  \|\Delta_{R^{(k)}} f\|_{L^2(\mu)}.
\end{align*}
Therefore, using Lemma \ref{lem:es1} and Lemma \ref{lem:es2} we see that
\begin{align*}
|\theta_t(1_{R^{(k)} \setminus R^{(k-1)}} \Delta_{R^{(k)}}f)(x)| + |&\theta_t( 1_{(R^{(k-1)})^c}B_{R^{(k-1)}}  b_{(R^{(k)})^a})(x)| \\
&\lesssim 2^{-\alpha k/2} \mu(R^{(k-1)})^{-1/2} \|\Delta_{R^{(k)}} f\|_{L^2(\mu)}, \qquad (x,t) \in W_R.
\end{align*}

Next, we note that
\begin{align*}
&\sum_{R:\, \ell(R) < 2^{s-r}} \mu(R) \Big[ \sum_{k=r+1}^{s+\textup{gen}(R)} 2^{-\alpha k /2} \mu(R^{(k-1)})^{-1/2} \|\Delta_{R^{(k)}} f\|_{L^2(\mu)} \Big]^2 \\
&\lesssim \sum_{R:\, \ell(R) < 2^{s-r}}  \mu(R)   \sum_{k=r+1}^{s+\textup{gen}(R)} 2^{-\alpha k /2} \mu(R^{(k-1)})^{-1} \|\Delta_{R^{(k)}} f\|_{L^2(\mu)}^2 \\
&= \sum_{k=r+1}^{\infty} 2^{-\alpha k /2} \sum_{m=k-s}^{\infty} \sum_{S:\, \ell(S) = 2^{k-m-1}} \|\Delta_{S^{(1)}} f\|_{L^2(\mu)}^2 \mu(S)^{-1} \mathop{\sum_{R:\,\ell(R) = 2^{-m}}}_{R \subset S} \mu(R) \\
&= \sum_{k=r+1}^{\infty} 2^{-\alpha k /2} \sum_{m=k-s}^{\infty} \sum_{S:\, \ell(S) = 2^{k-m-1}} \|\Delta_{S^{(1)}} f\|_{L^2(\mu)}^2 \\
&\lesssim  \sum_{k=r+1}^{\infty} 2^{-\alpha k /2} \sum_{m=k-s}^{\infty} \sum_{S:\, \ell(S) = 2^{k-m}} \|\Delta_S f\|_{L^2(\mu)}^2 \lesssim \sum_{S:\, \ell(S) \le 2^s} \|\Delta_S f\|_{L^2(\mu)}^2  \lesssim \|f\|_{L^2(\mu)}^2.
\end{align*}
Therefore, the first two terms of the splitting \eqref{eq:split1} are in control. The last part will be included in a paraproduct to be dealt with later.

\subsection{The case $(R^{(k-1)})^a = R^{(k-1)}$} We begin by writing
\begin{align*}
\Delta_{R^{(k)}} f = 1_{R^{(k)} \setminus R^{(k-1)}} \Delta_{R^{(k)}}f + 1_{R^{(k-1)}} \Delta_{R^{(k)}}f.
\end{align*}
The first term is in check by the argument above. We then decompose
\begin{align*}
 1_{R^{(k-1)}} \Delta_{R^{(k)}}f = \frac{\langle f \rangle_{R^{(k-1)}}}{\langle b_{R^{(k-1)}}\rangle_{R^{(k-1)}}}&b_{R^{(k-1)}} - \frac{\langle f \rangle_{R^{(k)}} }{\langle b_{(R^{(k)})^a}\rangle_{R^{(k)} }}b_{(R^{(k)})^a} \\
 &+ 1_{(R^{(k-1)})^c} \frac{\langle f \rangle_{R^{(k)}} }{\langle b_{(R^{(k)})^a}\rangle_{R^{(k)} }}b_{(R^{(k)})^a}.
\end{align*}
For the last term we have from Lemma \ref{lem:es2} that
\begin{displaymath}
|\theta_t( 1_{(R^{(k-1)})^c} b_{(R^{(k)})^a})(x)| \lesssim 2^{-\alpha k /2}, \qquad (x,t) \in W_R.
\end{displaymath}
For the term in front we simply use:
\begin{displaymath}
\frac{|\langle f \rangle_{R^{(k)}}| }{|\langle b_{(R^{(k)})^a}\rangle_{R^{(k)} }|} \lesssim |\langle f \rangle_{R^{(k)}}|.
\end{displaymath}
To finish the estimation of this term we bound
\begin{align*}
\sum_{R:\, \ell(R) < 2^{s-r}} \mu(R) \Big[& \mathop{\sum_{k=r+1}^{s+\textup{gen}(R)}}_{(R^{(k-1)})^a = R^{(k-1)}} 2^{-\alpha k /2} |\langle f \rangle_{R^{(k)}}| \Big]^2 \\
& \lesssim \sum_{S: \, \ell(S) \le 2^s} A_S |\langle f \rangle_S|^2 \lesssim \|f\|_{L^2(\mu)}^2.
\end{align*}
The last bound follows since \begin{displaymath}
A_S := \mathop{\sum_{S' \in \textup{ch}(S)}}_{(S')^a = S'} \mu(S')
\end{displaymath}
is a Carleson sequence. The rest will again become part of the paraproduct, which we will deal with in the next subsection.

\subsection{The bound for the paraproduct}
Combining the above two cases and collapsing the remaining telescoping summation we have reduced to estimating 
\begin{align*}
\mathop{\mathop{\sum_{R \in \mathcal{D}_{\textup{good}}}}_{\ell(R) < 2^{s-r}}}_{R \subset Q^*}& \iint_{W_R} \Big| \frac{\langle f \rangle_{R^{(r)}}}{\langle b_{(R^{(r)})^a} \rangle_{R^{(r)}}} \theta_t b_{(R^{(r)})^a}(x)
\Big|^2\,d\mu(x)\frac{dt}{t} \\
&\lesssim \sum_S |\langle f \rangle_S|^2 \mathop{\sum_{R \in \mathcal{D}}}_{S = R^{(r)}}\iint_{W_R} |\theta_t b_{S^a}(x)|^2 d\mu(x)\frac{dt}{t}.
\end{align*}
The key simplification is that this will not be handled using a Carleson estimate (which would be uniform i.e. we would be essentially bounding the above term for all functions). Instead, we rely on
the fact that the stopping time includes the Calder\'on--Zygmund stopping data of this specific function $f$:
\begin{align*}
 \sum_S &|\langle f \rangle_S|^2 \mathop{\sum_{R \in \mathcal{D}}}_{S = R^{(r)}}\iint_{W_R} |\theta_t b_{S^a}(x)|^2 d\mu(x)\frac{dt}{t}\\ 
&= \sum_{F \in \mathcal{F}} \sum_{S:\, S^a = F}  |\langle f \rangle_S|^2 \mathop{\sum_{R \in \mathcal{D}}}_{S = R^{(r)}}\iint_{W_R} |\theta_t b_{F}(x)|^2 d\mu(x)\frac{dt}{t} \\
&\lesssim \sum_{F \in \mathcal{F}} \alpha_f(F)^2 \sum_{S:\, S^a = F}  \mathop{\sum_{R \in \mathcal{D}}}_{S = R^{(r)}}\iint_{W_R} |\theta_t b_{F}(x)|^2 d\mu(x)\frac{dt}{t} \\
& \le \sum_{F \in \mathcal{F}} \alpha_f(F)^2 \sum_{R:\, R \subset F} \iint_{W_R} |\theta_t b_{F}(x)|^2 d\mu(x)\frac{dt}{t} \\
&\le \sum_{F \in \mathcal{F}} \alpha_f(F)^2 \iint_{\widehat F} |\theta_t b_{F}(x)|^2 d\mu(x)\frac{dt}{t}  \\
&\lesssim \sum_{F \in \mathcal{F}} \alpha_f(F)^2\mu(F) \lesssim \|f\|_{L^2(\mu)}^2.
\end{align*}
This estimate completes the proof of our main theorem, Theorem \ref{thm:main}.

\end{document}